\documentclass[12pt]{article}
\usepackage{bez123,calc,curves,ebezier,epic,eepic,graphicx,multiply,rotating}
\usepackage{algorithmic}
\everymath{\displaystyle}
\usepackage{graphicx}
\usepackage{pst-all}
\usepackage{color}
\usepackage{amssymb}
\usepackage{amsmath}
\newtheorem{prethm}{{\bf Theorem}}

\newenvironment{thm}{\begin{prethm}{\hspace{-0.5
               em}{\bf .}}}{\end{prethm}}
\newtheorem{prelemma}{{\bf Lemma}}

\newtheorem{preex}{{\bf Example}}

\newtheorem{preprop}{{\bf Proposition}}

\newenvironment{prop}{\begin{preprop}{\hspace{-0.5em}{\bf .}}}{\end{preprop}}
\newtheorem{precor}{{\bf Corollary}}

\newenvironment{cor}{\begin{precor}{\hspace{-0.5
               em}{\bf .}}}{\end{precor}}
\newtheorem{preremark}{{\bf Remark}}

\newenvironment{remark}{\begin{preremark}{\hspace{-0.5
               em}{\bf.}}}{\end{preremark}}
\newtheorem{preprob}{{\bf Problem}}

\newtheorem{predefin}{{\bf Definition}}

\newtheorem{preconj}{{\bf Conjecture}}

\newtheorem{preprobb}{{\bf Problem}}

\newtheorem{prelem}{{\bf Theorem}}

\newenvironment{proof}{{\bf Proof.}\rm }{\hfill{$\Box$}}

\newtheorem{presolution}{{\bf Solution.}}

\def\newpic#1{}
\def\qed{\ifhmode\unskip\nobreak\fi\quad\ifmmode\Box\else$\Box$\fi}

\title{\Large\bf\noindent On dynamic monopolies of graphs:
the average and strict majority thresholds}

\author{\large\bf Kaveh Khoshkhah, Hossein Soltani, Manouchehr Zaker\footnote{Corresponding author: mzaker@iasbs.ac.ir}
\vspace{5mm}\\
    Department of Mathematics,\\
     Institute for Advanced Studies in Basic Sciences,\\
    Zanjan 45137-66731, Iran}
    \date{}
\begin{document}
\maketitle
\begin{abstract}
\noindent Let $G$ be a graph and ${\mathcal{\tau}}: V(G)\rightarrow \Bbb{N}\cup \{0\}$ be an assignment of thresholds to the vertices of $G$. A subset of vertices $D$ is said to be a dynamic monopoly corresponding to $(G, \tau)$ if the vertices of $G$ can be partitioned into subsets $D_0, D_1, \ldots, D_k$ such that $D_0=D$ and for any $i\in \{0, \ldots, k-1\}$, each vertex $v$ in $D_{i+1}$ has at least $\tau(v)$ neighbors in $D_0\cup \ldots \cup D_i$. Dynamic monopolies are in fact modeling the irreversible spread of influence in social networks. In this paper we first obtain a lower bound for the smallest size of any dynamic monopoly in terms of the average threshold and the order of graph. Also we obtain an upper bound in terms of the minimum vertex cover of graphs. Then we derive the upper bound $|G|/2$ for the smallest size of any dynamic monopoly when the graph $G$ contains at least one odd vertex, where the threshold of any vertex $v$ is set as $\lceil (deg(v)+1)/2 \rceil$ (i.e. strict majority threshold). This bound improves the best known bound for strict majority threshold. We show that the latter bound can be achieved by a polynomial time algorithm. We also show that $\alpha'(G)+1$ is an upper bound for the size of strict majority dynamic monopoly, where $\alpha'(G)$ stands for the matching number of $G$. Finally, we obtain a basic upper bound for the smallest size of any dynamic monopoly, in terms of the average threshold and vertex degrees. Using this bound we derive some other upper bounds.
\end{abstract}




\section{Introduction}

\noindent All graphs in this paper are undirected graphs without multiple edges or loops. For any graph $G$ we denote the vertex set, the edge set and the order of $G$ by, $V(G)$, $E(G)$ and $|G|$, respectively. For other graph theoretical notations not defined in this paper we refer the reader to \cite{BM}. In this paper we denote the set of positive natural numbers by $\Bbb{N}$. Let $G$ be a graph and ${\mathcal{\tau}}: V(G)\rightarrow \Bbb{N}\cup \{0\}$ be an assignment of thresholds to the vertices of $G$ such that $\tau(v)\leq deg(v)$, where $deg(v)$ is the degree of $v$ in $G$. A subset of vertices $D$ is said to be a {\it $\tau$-dynamic monopoly} or simply {\it dynamic monopoly} if the vertices of $G$ can be partitioned into subsets $D_0, D_1, \ldots, D_k$ such that $D_0=D$ and for any $i$, $1\leq i \leq k-1$, each vertex $v$ in $D_{i+1}$ has at least $\tau(v)$ neighbors in $D_0\cup \ldots \cup D_i$. The usual formulation of dynamic monopolies is in terms of a discrete time dynamic process defined as follows. Consider a dynamic process on the vertices of $G$, where some vertices of $G$ are considered as active vertices at the beginning of the process. Denote the set of active vertices at any discrete time $t\geq 0$ by $D_t$. Assume that at the beginning of the process (i.e. at time zero), the vertices of a subset $D\subseteq V(G)$
are active. Hence $D_0=D$. At each discrete time $i$ any inactive vertex $v$ is activated provided that $v$ has at least $\tau(v)$ active neighbors in $D_0 \cup \ldots \cup D_{i-1}$. If at the end of the process all vertices are active then the starting subset $D$, is called dynamic monopoly or simply {\it dynamo}. Given $(G, \tau)$, by the average threshold in $G$ we mean $\sum \tau(v)/|G|$. Some well-known threshold assignments for the vertices of a graph $G$ are simple and strict majority thresholds. In {\it simple majority threshold} we set $t(v)=deg(v)/2$ for any vertex $v$ of $G$ and in {\it strict majority threshold} we have $t(v)=\lceil (deg(v)+1)/2 \rceil$. In this paper by {\it a strict majority dynamic monopoly} we mean any dynamo corresponding to the strict majority threshold. In the recent years there has been a great interest to study of the dynamic monopolies in graphs. Strict majority dynamos in some special families of graphs were studied in \cite{FGS, FKRRS, FLLPS, LPS}. The first complexity results concerning dynamic monopolies with general thresholds appeared in \cite{C}. The first theoretical results for graphs with general thresholds were obtained in \cite{Z}. Dynamic monopolies with constant thresholds were studied in \cite{DR}, where some hardness and algorithmic results have also been obtained. An important motivation to study the dynamic monopolies is their applications in formulation of the spread of influence in social networks \cite{C, KKT}. Some examples of these phenomena are
the spread of virus among a population or in a web of computers, spread of innovation
or a new product in a community, spread of opinion in elections and etc. Dynamic monopolies have also applications in viral marketing \cite{DRi}. Dynamic monopolies of random graphs were studied in \cite{CL2, K}. In \cite{CL}, the authors studied the dynamic monopolies with strict majority thresholds in undirected and directed graphs. Dynamic monopolies of graph products were studied in \cite{ABST}. More studies on dynamic monopolies can be found in \cite{B, P2}.

\noindent {\bf The outline of the paper:} In Section 2 we introduce some concepts involving the average threshold and prove some basic results and bounds for the size of dynamic monopolies with given average thresholds. Section 3 devotes to the study of strict majority dynamic monopolies. We first derive the upper bound $|G|/2$ for the smallest size of any strict majority dynamic monopoly when the graph $G$ contains at least one odd vertex. This bound improves the best known bound for strict majority threshold. We show that the latter bound can be achieved by a polynomial time algorithm. Also in Section 3 we show that $\alpha'(G)+c$ is an upper bound for the size of strict majority dynamic monopoly of any graph $G$ with $c$ connected components, where $\alpha'(G)$ stands for the matching number of $G$. In Section 4, we obtain an upper bound for the smallest size of any dynamic monopoly, in terms of the average threshold and vertex degrees. Using this bound we show that given any graph $(G, \tau)$ on $n$ vertices and with average threshold $\bar{t}$, then there exists a $\tau$-dynamic monopoly with at most $\frac{n\bar{t}}{\delta(G)+1}$ vertices, where $\delta(G)$ is the minimum degree of $G$. We show that this bound is achieved by a polynomial time algorithm.

\section{Average thresholds}

\noindent By a threshold assignment to the vertices of a graph $G$ we mean any function $\tau: V(G) \rightarrow \Bbb{N}\cup \{0\}$ such that the threshold of any vertex $v$ is at most $deg(v)$, where $deg(v)$ stands for the degree of $v$ in $G$. We denote the threshold of a vertex $v$ by $t(v)$. Let a graph $G$, a threshold assignment $\tau$ and a subset $M\subseteq V(G)$ be given. For any $i$, $i=0, 1, \ldots$, we define a subset $D_i$ as follows. Set $D_0=M$. Assume that $D_0, D_1, \ldots, D_i$ have been defined. Define $D_{i+1}$ as a subset consisting of all vertices $v$ such that $v$ has at least $t(v)$ neighbors in $D_0\cup D_1 \cup \ldots \cup D_i$. It is possible that $D_i=\varnothing$ for some $i\geq 1$. By the {\it $\tau$-dynamic process} starting from $M$ we mean the sequence $D_0, D_1, \ldots$ If the sequence is such that $V(G) = \bigcup_i D_i$, then $M$ is called a $\tau$-dynamic monopoly (or simply dynamic monopoly). Given a graph $G$ and a threshold assignment $\tau$ for its vertices, we denote the minimum number of vertices in any $\tau$-dynamic monopoly of $G$ by $dyn_{\tau}(G)$. Denote the maximum (resp. minimum) threshold in $G$ by $t_M$ (resp. $t_m$). In \cite{Z}, some bounds in terms of the minimum or maximum thresholds for the smallest size of dynamic monopolies in graphs
were obtained. Also in \cite{Z}, graphs with probabilistic thresholds were considered and the importance of the expectation (or average) of thresholds in lower-bounding the size of dynamic monopolies was shown. It is more useful to obtain bounds in terms of the average threshold. For any threshold assignment $\tau$ of a graph $G$, by the average
threshold of $\tau$ we mean $\sum_{v\in G} \tau(v)/|G|$ and denote it by $\overline{\tau}$. In applications too the average threshold is more accessible than the minimum or maximum thresholds. In other words, in most applications our knowledge is only about the average of thresholds in a network. In this section we
intend to study the dynamic monopolies and extend some previous results in terms of the average threshold. For any rational number $t$ we introduce $Dyn_{_{\bar{t}=t}}(G)$ which is the main parameter to be studied in this paper and is defined as follows, where the maximum is taken over all threshold assignments $\tau$ such that $\overline{\tau}=t$:

$$Dyn_{_{\bar{t}=t}}(G)=\max_{\tau: \overline{\tau}=t}~ dyn_{\tau}(G).$$

\noindent Assume that a family ${\mathcal{F}}$ of graphs is given such that any graph $G$ from ${\mathcal{F}}$ is equipped with a threshold assignment. Recall that the smallest size of any dynamo of $G$ corresponding to its threshold assignment is denoted by $dyn(G)$. In \cite{Z}, the family ${\mathcal{F}}$ is called dynamo-unbounded if there exists a function $f(x)$ satisfying $f(x)\rightarrow \infty$ as $x\rightarrow \infty$ such that for any graph $G$ from $\mathcal{F}$ one has $f(n)\leq dyn(G)$, where $n=|G|$. The following result was proved in \cite{Z}. Let $(G, \tau)$ be a graph of order $n$ and $\epsilon(G)$ be the edge density of $G$, i.e. $\epsilon(G)=|E(G)|/|G|$. Set $t=\min \{t(v):
v\in V(G)\}$. Then $n(1-\frac{\epsilon(G)}{t})\leq dyn(G)$. Using the latter bound, some families were proved to be dynamo-unbounded in \cite{Z}. In this section our aim is to generalize the latter lower bound in terms of average threshold.

\begin{thm}
Let $G$ be a graph with maximum degree $\Delta$ and $\tau$ be a threshold assignment to the vertices of $G$.
Let also $\bar{t}$ and $t_M$ denote the average and maximum threshold of $\tau$, respectively. For any $\tau$-dynamic monopoly $M$ of $G$
we have $$|M|\geq |G|(1-\frac{\epsilon(G)}{\bar{t}})(\frac{\bar{t}}{t_M}) \geq |G|(1-\frac{\epsilon(G)}{\bar{t}})(\frac{\bar{t}}{\Delta}).$$\label{ave}
\end{thm}

\noindent \begin{proof} Let $n=|G|$. There exists a partition $D_0\cup D_1 \cup \ldots \cup D_t$ of $V(G)$ such that $D_0=M$ and for any $i\geq 1$ and any vertex $v\in D_i$, there exist at least $t(v)$ edges between $v$ and $D_0\cup \ldots \cup D_{i-1}$. Therefore at least $\sum_{v\in V(G)\setminus M} t(v)$ edges exist in $G$. We have now $\sum_{v\in V(G)} t(v) - t_M|M| \leq \sum_{v\in V(G)\setminus M} t(v) \leq |E(G)|$. It follows that $\bar{t} - (t_M/n) |M| \leq \epsilon(G)$ or
$n/t_M(\bar{t}-\epsilon(G))\leq |M|$. This completes the proof of the first inequality.
The second one is easily obtained by $t_M \leq \Delta$.
\end{proof}

\noindent Let ${\mathcal{F}}$ be any family of graphs such that for some positive constant $\delta$, $\min\{t(v): v \in G\} \geq \epsilon(G) + \delta$ for any graph $G\in {\mathcal{F}}$. Then as shown in \cite{Z}, ${\mathcal{F}}$ is dynamo-unbounded family. The following corollary gives a lower bound in terms of the edge density of graphs.

\begin{cor}
Let $\delta$ be any positive constant and $G$ any graph with edge density $\epsilon$. Let also $\bar{t}$ be any constant with $\bar{t}\geq (1+\delta)\epsilon$. Let $\tau$ be any threshold assignment with average $\bar{t}$ and $M$ be any $\tau$-dynamic monopoly for $G$. Then
$$|M|\geq \delta\epsilon.$$
\end{cor}

\noindent \begin{proof}
The proof is easily obtained using the lower bound of Theorem \ref{ave} and that
$\bar{t}\geq (1+\delta)\epsilon$ and $t_M < n$.
\end{proof}

\noindent It is worth-mentioning that when $\epsilon / \bar{t}\rightarrow 1$ in a family of graphs i.e. when the lower bound of Theorem \ref{ave}
tends to zero then it is possible that the dynamic monopoly of all members of that family is bounded by a constant number. As a simple example, consider the family of complete graphs $K_{n}$ where $n=1, 2, \ldots$, and let the thresholds in $K_n$ be $1, 2, 3, \ldots, n-1, n-1$. Note that $\epsilon(K_n)=(n-1)/2$ and the average threshold is $\epsilon(K_n)+ (n-1)/n$. It follows that $\epsilon(K_n)/\bar{\tau}(K_n)\rightarrow 1$ as $n\rightarrow \infty$. From other side, it is easy to see that a single vertex with threshold $n-1$ in $K_n$ forms a dynamo for $K_n$.

\noindent Let $G$ be a $(2r+1)$-regular graph on $n$ vertices and $t(v)=r+1$, for every vertex $v$ of $G$. Then as shown in \cite{Z}
any dynamo for $G$ has at least $(n+2r)/(2r+2)$ vertices. Using Theorem \ref{ave}
we have the following corollary.

\begin{cor}
Let $G$ be a $(2r+1)$-regular graph on $n$ vertices. Let also $\tau$ be a threshold assignment for $G$ with average threshold $r+1$. Then any dynamo for $(G,\tau)$ has at least $n/(4r+2)$ vertices.
\end{cor}

\noindent In the following proposition we obtain a general upper bound for $Dyn_{_{\bar{t}=2\epsilon}}(G)$ when the average threshold ${\bar{t}}$ is any arbitrary value such that $0\leq \bar{t} \leq 2\epsilon(G)$. Note that the average threshold in a graph $G$ can not exceed $2\epsilon(G)$. In the following we denote the smallest size of any vertex cover of $G$ by $\beta(G)$. Recall that a vertex cover is a subset $S$ of vertices such that any edge of the graph has at least one endpoint in $S$.

\begin{prop}
Let $G$ be a graph without isolated vertices.

\noindent (i) Let $\tau$ be any threshold assignment with average $2\epsilon$. Then any $\tau$-dynamo has $\beta(G)$ vertices. In particular $Dyn_{_{\bar{t}=2\epsilon}}(G)=\beta(G)$.

\noindent (ii) For any constant $t\leq 2\epsilon$, $Dyn_{_{\bar{t}=t}}(G)\leq \beta(G)$.\label{beta}

\end{prop}

\noindent \begin{proof}
Let $t_1, t_2, \ldots, t_n$ be any set of thresholds such that their average is $2\epsilon$. Let also $M$ be any dynamic monopoly for these thresholds. Since $t_i\leq deg(v_i)$ and the
average of thresholds is $2\epsilon$ then $t_i=deg(v_i)$ for any $i$. It is now clear that from any edge $e=uv$ of the graph either $u$ or $v$ should be in $M$, since otherwise neither $u$ nor $v$ will become active until the end of the process. Hence $M$ is a vertex cover. From other side any vertex cover is a dynamo. Therefore
the cardinality of $M$ should be $\beta$. This proves part (i).

\noindent Now we prove part (ii). Assume that $Dyn_{_{\bar{t}=t}}(G)$ is achieved
by a specific set of thresholds $t_1, \ldots, t_n$. Let $K$ be a vertex cover for $G$.
Note that $K$ is a dynamo for the thresholds $t_1, \ldots, t_n$. By definition
$Dyn_{_{\bar{t}=t}}(G)$ is the size of smallest dynamo in $G$. Hence $Dyn_{_{\bar{t}=t}}(G)\leq |K| = \beta$.
\end{proof}

\noindent In a graph $G$ denote the maximum number of independent vertices and the chromatic number of $G$ by $\alpha(G)$ and $\chi(G)$, respectively. It is a well-known fact (see e.g. \cite{BM}) that $\alpha(G)+\beta(G)=|G|$. Also it is easily seen that $|G|\leq \alpha(G)\chi(G)$. By these notes and Proposition \ref{beta}, the following corollary is easily obtained.

\begin{cor}
For any constant $t\leq 2\epsilon$, $Dyn_{_{\bar{t}=t}}(G)\leq |G|(1-\frac{1}{\chi(G)})$.\label{chi}
\end{cor}

\noindent Corollary \ref{chi} shows that if we consider the family of graphs with bounded maximum degree then for some constant $\lambda$, $Dyn_{\bar{t}}(G)\leq \lambda |G|$ for any graph $G$ from the family. The following proposition shows that when $\Delta\rightarrow \infty$ then there exists no upper bound for $Dyn_{\bar{t}}(G)$ in the form of $\lambda |G|$, where $\lambda$ is a constant strictly less than one.

\begin{prop}
There exists an infinite sequence of graphs $G_1, G_2, \ldots$ such that $|G_n|\rightarrow \infty$ and $$\lim_{n\rightarrow \infty} \frac{Dyn_{\epsilon}(G_n)}{|G_n|} =1.$$\label{example}
\end{prop}

\noindent \begin{proof}
Set $G_1=K_2$. For any integer $n\geq 2$ we construct a graph denoted by $G_n$ on $n(n-1)+n[n(n-1)-1]$ vertices as follows. We consider one copy of $K_{n(n-1)}$ and $n(n-1)-1$ vertex disjoint copies of $K_n$. There exists no edge between these copies of $K_n$ but we connect any vertex from any copy of $K_n$ to any vertex of $K_{n(n-1)}$ by an edge. Therefore $G_n$ contains ${n(n-1)\choose 2}+ [n(n-1)-1]{n\choose 2}+n(n-1)[n(n-1)-1]n$ edges. We simplify the latter value and obtain $|E(G_n)|=(n^3-n)[n^2-n-1]$. Now we obtain a suitable threshold assignment for $G_n$ with average $\epsilon(G_n)$. For any vertex $v$ from the copy $K_{n(n-1)}$ of $G_n$ set $t(v)=0$. For any vertex $u$ from the $K_n$ copies of $G_n$ set $t(u)=deg_{G_n}(u)=[n(n-1)+(n-1)]$. We obtain that $\sum_{v\in G_n}t(v)=|E(G_n)|$ and therefore the average threshold is exactly $\epsilon(G_n)$.

\noindent Now let $M$ be any dynamic monopoly corresponding to the given threshold assignment of $G_n$. Noting that the threshold of any vertex in each copy of $K_n$ in $G_n$ is its degree in the whole graph then we obtain that $M$ should contain a vertex cover from each $K_n$ copy of $G_n$. It implies that $M\geq [n(n-1)-1](n-1)$. Hence we have $Dyn_{\epsilon}(G_n)\geq [n(n-1)-1](n-1)$. Finally
$$1\geq \lim_{n\rightarrow \infty} \frac{Dyn_{\epsilon}(G_n)}{|G_n|}\geq \lim_{n\rightarrow \infty}~\frac{[n(n-1)-1](n-1)}{n(n-1)+n[n(n-1)-1]}=1.$$
\noindent This completes the proof.
\end{proof}\\

\section{Strict majority dynamic monopolies}

\noindent In this section we consider graphs with strict majority thresholds, i.e.
for any vertex $v$ we set $t(v)=\lceil (deg(v)+1)/2 \rceil$. In \cite{CL} it was shown that any graph $G$ contains a strict majority dynamo of at most $\lceil |G|/2 \rceil$ vertices. In Corollary \ref{bound-major} we improve their result. Strict majority dynamic monopolies were also studied in \cite{ABW}, where the same bound as in \cite{CL} were presented.

\noindent Let $G$ be a graph and $\sigma$ any vertex ordering of $G$. Denote the order of a vertex $v$ in $G$ by $\sigma(v)$. For any two vertices $u$ and $v$, $\sigma(u)<\sigma(v)$
means that $u$ appeares before $v$ in the ordering $\sigma$. Also denote the neighborhood set of any vertex $v$ by $N(v)$. For any vertex $v$ we define $f_{\sigma}(v)$ as follows
$$f_{\sigma}(v)=|N(v)\cap \{u: \sigma(u)>\sigma(v)\}|-|N(v)\cap \{u: \sigma(u)<\sigma(v)\}|$$

\begin{thm}
\noindent Let $G$ be any connected graph.

\noindent (i) If $G$ contains at least one vertex of odd degreed, then there exists an ordering $\sigma$ such that for any vertex $v$, $f_{\sigma}(v)\not= 0$.

\noindent (ii) If all degrees in $G$ are even, then there exists an ordering $\sigma$ such that
$f_{\sigma}(v)\not= 0$ for all but at most one vertex $v$ of $G$. Moreover, in the case that for some vertex $v$, $f_{\sigma}(v)= 0$ then $v$ can be taken as any arbitrary vertex of $G$.\label{order}
\end{thm}

\noindent \begin{proof}
We prove the following stronger claim:

\noindent {\bf Claim:} There exists an ordering $\sigma$ satisfying the conditions of the theorem which has also the following stronger property. For any $u$, $v$ and $w$ if $f_{\sigma}(u)>0$, $f_{\sigma}(v)<0$ and $f_{\sigma}(w)=0$ then $\sigma(u)< \sigma(w)<
\sigma(v)$.

\noindent We prove the claim by induction on $|G|$. The assertion trivially holds
when $|G|$ is 1 or 2. Assume that it holds for all graphs of less than $n$ vertices
and let $G$ be a graph with $|G|=n$. If $G$ contains a vertex of odd degree then we let $x$ be a vertex of degree odd in $G$, otherwise let $x$ be an arbitrary vertex of $G$. Let $A_1, \ldots, A_k$ be the connected components of $G\setminus x$. By the induction hypothesis for each $A_i$ there corresponds an ordering $\sigma^i$ such that the associated function $f_{\sigma^i}$ satisfies the conditions of the claim and if there exists a vertex say $u$ in $A_i$ whose $f_{\sigma^i}$ is zero then $u$ can be chosen as a neighbor of $x$ in $G$ (since the vertex with $f_{\sigma}=0$ can be taken as any arbitrary vertex in $G$). Let $A_i^+$ (resp. $A_i^-$) be the sequence of vertices in $A_i$ whose $f_{\sigma^i}$ is positive (resp. negative). Let also $A_i^0$ be the vertex (if exists) in $A_i$ with
$f_{\sigma^i}=0$. Now we define an ordering $\sigma$ on $V(G)$ obtained by the sequence of
vertices specified in the following list from left to right
$$A_1^+, \ldots, A_k^+, A_1^0, \ldots, A_k^0, x, A_1^-, \ldots, A_k^-.$$
\noindent We note that the order of vertices in $A_i$ in both orders $\sigma$ and $\sigma^i$ is the same. Let $u$ be any vertex of $A_i$. Aside from $A_i$ itself, $u$ can only be adjacent to $x$. Since the position of $x$ in $\sigma$ is after $A_i^+$ and before $A_i^-$ then the sign of $f_{\sigma}(u)$ is the same as sign of $f_{\sigma^i}(u)$. Assume that there exists $u\in A_i$ with $f_{\sigma^i}(u)=0$. Since in $\sigma$, $x$ is appeared after $u$ then $f_{\sigma}(u)\not= 0$.
\end{proof}

\noindent As we mentioned before, it was shown in \cite{CL} that $G$ contains a strict majority dynamic monopoly
of cardinality at most $\lceil |G|/2 \rceil$. The following corollary gives a stronger result that if $G$ has at least one vertex of odd degree, than it admits a dynamo of at most
$n/2$ vertices. Also the proof is simpler than that of \cite{CL}.

\begin{cor}
Let $G$ be a graph on $n$ vertices and $\tau$ a threshold function obtained by $t(v)=\lceil (deg(v)+1)/2 \rceil$
for any vertex $v$. Then there exists a $\tau$-dynamic monopoly $M$ such that $|M|\leq \lceil n/2 \rceil$. Moreover, if $G$ contains a vertex of degree odd then there exists such a set $M$ with $|M|\leq n/2$.\label{bound-major}
\end{cor}

\noindent \begin{proof}
Let $\sigma$ be an ordering of the vertices of $G$ satisfying the conditions of Theorem \ref{order}. Let $M$ be the set of vertices $v$ such that $f(v)\geq 0$. We observe that $M$ is a strict majority dynamo. In fact the vertices with negative $f$ become active in turn according to their order in $\sigma$. Similarly if $M$ is the set of vertices $v$ with $f(v)\leq 0$ then $M$ is a dynamo. The vertices with positive $f$ become active in turn according to reverse of their order in $\sigma$. Now at least one of these sets have no more that $\lceil n/2 \rceil$ vertices.
\end{proof}

\noindent The following remark is immediate from Theorem \ref{order} and the proof of Corollary \ref{bound-major}.

\begin{remark}
Let $G$ be a graph on even number of vertices and $v$ be any vertex of $G$. Then $G$ admits a strict majority dynamo with cardinality at most $|G|/2$ which contains the vertex $v$.\label{remar}
\end{remark}

\noindent The methodology of the proof of Theorem \ref{order} shows that there exists a polynomial time recursive algorithm which constructs the ordering satisfying the conditions of Theorem \ref{order}. Using this ordering Corollary \ref{bound-major} easily obtains a dynamic monopoly with at most $\lceil n/2 \rceil$ vertices for any graph of order $n$. We have therefore the following remark.

\begin{remark}
There exists a polynomial time algorithm which for any connected graph $G$ on $n$ vertices, outputs a strict majority dynamo with at most $\lceil n/2 \rceil$ vertices.
\end{remark}

\noindent We end this section with relating the majority strict dynamic monopolies to matching number of graphs. By the matching number of $G$ we mean the maximum number of independent edges in $G$. In obtaining the next result we shall make use of a theorem from \cite{KSNZ}. For this purpose we need some terminology. By a graph parameter $p$ we mean any function $p$ from the set of all graphs to non-negative integers such that if $G$ and $H$ are two isomorphic graphs then $p(G)=p(H)$. Also a graph parameter $p$ is called subadditive if $p(G\cup H)\leq p(G)+p(H)$, where $G\cup H$ is the vertex disjoint union of two graphs $G$ and $H$. The following was proved in \cite{KSNZ}.

\begin{thm}
Let $p$ be any subadditive graph parameter such that for any graph $G$ and any vertex $v\in G$, $p(G)\leq p(G\setminus v)+1$. Assume that there exists a constant $1\leq t <2$ such that for any graph $G$ on odd number of vertices, $p(G)\leq t(|G|-1)/2$. Then $p(G)\leq \lfloor t \alpha'(G) \rfloor$ for any graph $G$.\label{parameter}
\end{thm}

\noindent Our result is as follows, where by $\alpha'(G)$ we mean the maximum number of independent edges.

\begin{thm}
Any connected graph $G$ admits a strict majority dynamo of size at most $\alpha'(G)+1$. Furthermore if $G$ has $c$ connected components then $G$ admits a strict majority dynamo of size at most $\alpha'(G)+c$.
\end{thm}

\noindent \begin{proof}
We define a graph parameter $p$ as follows. For any connected graph $G$ we define:
\begin{equation*}
p(G) =
\begin{cases}
dyn(G)-1 & \text{if } dyn(G)= \lceil \frac{|G|+1}{2} \rceil,\\
dyn(G) & \text{if } otherwise.
\end{cases}
\end{equation*} where $dyn(G)$ is the smallest size of any strict majority dynamo in $G$.
\noindent For a non-connected graph $G$ consisting of the connected components $H_1, H_2, \ldots, H_k$ we define $p(G)=\sum_{i}p(H_i)$. Note that $p$ is indeed a subadditive graph parameter. Note also that $p(K_1)=0$ since $dyn(K_1)=1$. In the following we show that $p$ satisfies the Lipschitz inequality: $p(G)\leq p(G\setminus v)+1$, for any $v\in G$. It is easily seen by subadditivity of $p$ and $(G \cup H)\setminus v=(G\setminus v) \cup H$ for any $v\in G$ that it is enough to prove the Lipschitz property for connected graphs.

\noindent Now let $G$ be a connected graph and $v\in G$. Assume that the connected components of $G\setminus v$ are $G_1, G_2, \ldots, G_t$. In the following we construct a strict majority dynamo $M$ for $G$ such that $v\in M$ and $|M\cap G_i| \leq |G_i|/2$. For this purpose, let $M$ be a strict majority dynamo containing $v$ with the smallest cardinality. We show that $|M\cap G_i| \leq |G_i|/2$, for any $i$. Assume by the contrary that for some $j$, $|M\cap G_j| > |G_j|/2$. Since $G_j$ itself admits a strict majority dynamo with at most $\lceil |G_j|/2 \rceil$ elements, then $|M\cap G_j|=\lceil (|G_j|+1)/2 \rceil$. Note that in this case $|G_j|$ is odd. We consider the subgraph of $G$ induced by $V(G_j)\cup \{v\}$. By Remark \ref{remar}, the latter graph admits a dynamo say $M_0$ containing $v$ with at most $\lceil (|G_j|+1)/2 \rceil$ elements. Note that $|M_0\setminus \{v\}|<|M\cap G_j|$. We obtain a new dynamo for $G$ as follows $M_{new}=(M\setminus G_j) \cup (M_0\setminus \{v\})$. Now $M_{new}$ is a dynamo containing $v$ and with cardinality less than $M$, a contradiction.

\noindent For any $i$, denote by $dyn'(G_i)$ the smallest size of any strict majority dynamic monopoly
in the subgraph of $G$ induced by $V(G_i)\cup \{v\}$, where the vertex $v$ is already an active vertex. Note that using the strict majority dynamo $M$ obtained in the previous paragraph we have $dyn'(G_i)\leq |G_i|/2$. In the following we show that $dyn'(G_i)\leq p(G_i)$. It is clear that $dyn'(G_i)\leq dyn(G_i)$. If $dyn'(G_i)= dyn(G_i)$ then since $dyn'(G_i)\leq |G_i|/2$,
$dyn(G_i)\leq |G_i|/2$ and so $p(G_i)=dyn(G_i)$ or $p(G_i)=dyn'(G_i)$. But if $dyn'(G_i)\leq dyn(G_i)-1$ then by the definition of $p(G_i)$, $dyn'(G_i)\leq p(G_i)$. We have now the following inequalities, where $M$ is the dynamo we obtained in the above paragraph

$$p(G)\leq dyn(G) \leq |M| = 1+\sum_i dyn'(G_i)\leq 1+\sum_i p(G_i) =p(G\setminus v)+1.$$

\noindent By Theorem \ref{parameter}, for any graph $G$ consisting of the connected components $G_1, G_2, \ldots, G_c$ we have $p(G)\leq \alpha'(G)$. From one side $p(G)=\sum_i p(G_i)$. From other side by the definition of $p$ for connected graphs, $p(G_i)\geq dyn(G_i)-1$. Combining these results with $p(G)\leq \alpha'(G)$ yields $dyn(G)\leq \alpha'(G)+c$.
\end{proof}

\section{Some upper bounds}

\noindent This section is devoted to presenting some upper bounds for the size of dynamic monopolies in terms of the average threshold. We first obtain a basic upper bound in terms of average threshold and vertex degrees. In the following theorem, for any vertex $v$ and subset $S$ of the vertices in a graph $G$, we denote the set of the neighbors of $v$ in $S$ by $N_S(v)$.

\begin{thm}
Let $G$ be a graph with degree sequence $d_1 \leq d_2 \leq \ldots \leq d_n$ in increasing
form. Let also $\tau$ be any threshold assignment for the vertices of $G$ with average threshold $\bar{t}$. Then $$Dyn_{_{\bar{t}=t}}(G)\leq \max \{k: \sum_{i=1}^k (d_i+1) \leq n\bar{t}\}.$$\label{upper}
\end{thm}

\noindent \begin{proof}
Denote the threshold and the degree of any vertex $v$ of $G$ by $t(v)$ and $deg(v)$, respectively. Let $M$ be any $\tau$-dynamic monopoly of $G$ with minimum cardinality. We partition $V(G)\setminus M$ into two subsets $A$ and $B=V(G)\setminus M\setminus A$ where
$A=\{v\in V(G)\setminus M: |N_M(v)|\leq t(v)\}$. Note that if $v\in B$ then $|N_M(v)|>t(v)$. We make the following claim:

\noindent {\bf Claim:} For any vertex $x\in M$, $|N_{A\cup B}(x)|\geq |N_B(x)|+deg(x)-t(x)+1$.

\noindent {\bf Proof of the claim:} Assume by the contrary that $|N_{A\cup B}(x)|\leq |N_B(x)|+deg(x)-t(x)$, for some vertex $x$. We have $|N_M(x)|=deg(x)-|N_{A\cup B}(x)|\geq deg(x)-|N_B(x)|+deg(x)-t(x)=t(x)-|N_B(x)|$. It follows that $|N_{M\cup B}(x)|\geq t(x)$.
This means that the vertex $x$ can be active by activation of all vertices in $M\cup B$. From the other side by the definition of the set $B$ all vertices of $B$ can be active by activation of $M\setminus \{x\}$. In other words $M\setminus \{x\}$ becomes a dynamo, which contradicts the minimality of $M$.

\noindent We have the following
$$\sum_{y\in A\cup B}|N_M(y)|=\sum_{y\in A}|N_M(y)|+\sum_{y\in B}|N_M(y)|\leq \sum_{y\in A}t(y)+\sum_{y\in B}|N_M(y)|,~~~~~{\bf (1)}$$
\noindent and from the claim
$$\sum_{x\in D}|N_{A\cup B}(x)|\geq \sum_{x\in D}|N_B(x)|+\sum_{x\in D}(deg(x)+1)-\sum_{x\in D}t(x).~~~~~~~~{\bf (2)}$$
\noindent We have the equality $\sum_{y\in A\cup B}|N_D(y)|=\sum_{x\in D}|N_{A\cup B}(x)|$.

\noindent Therefore from (1) and (2)
$$\sum_{x\in D}|N_B(x)|+\sum_{x\in D}(deg(x)+1)-\sum_{x\in D}t(x)\leq \sum_{y\in A}t(y)+\sum_{y\in B}|N_D(y)|.$$
\noindent Also the equality $\sum_{x\in D}|N_B(x)|=\sum_{y\in B}|N_D(y)|$ holds.
\noindent It follows that $$\sum_{x\in D}(deg(x)+1)-\sum_{x\in D}t(x)\leq \sum_{y\in A}t(y),$$
\noindent and finally $$\sum_{x\in D}(deg(x)+1)\leq \sum_{y\in A\cup D}t(y)\leq \sum_{v\in V(G)}t(v)=n\bar{t}.$$
\noindent We conclude that $|D|\leq \max \{k: \sum_{i=1}^k (d_i+1) \leq n\bar{t}\}$.
\end{proof}

\noindent In the following proposition we show that the bound of Theorem \ref{upper}
can be achieved by an efficient algorithm.

\begin{prop}
There exists an ${\mathcal{O}}(n^3)$ algorithm which for any graph $G$ on $n$ vertices
and any threshold assignment of $G$ with average $\bar{t}$, outputs a dynamo $M$ such that $|M|\leq \max \{k: \sum_{i=1}^k (d_i+1) \leq n\bar{t}\}$.
\end{prop}

\noindent \begin{proof} The description of the algorithm is as follows. At each time step of the algorithm we have a dynamic monopoly denoted by $M$. At the beginning, we set $M=V(G)$. We modify the set $M$ through the execution of the algorithm such that at each step of the algorithm, $M$ is a dynamo and at the last step we obtain a dynamo with the desired cardinality.

\noindent At each time we have a dynamic monopoly $M$ and its corresponding sets $A$ and $B$. The set $A$ is defined as $A=\{v\in V(G)\setminus M: |N_M(v)|\leq t(v)\}$. Also $B=V(G)\setminus M\setminus A$. Assume that at a certain step of the procedure we have
$M$, $A$ and $B$. We modify $M$ as follows. We scan all vertices in $M$ to find a vertex
$v$ such that $v$ does not satisfy the following condition
$$|N_{A\cup B}(v)|\geq |N_B(v)|+deg(v)-t(v)+1.$$
\noindent There are two possibilities:

\noindent {\bf Case 1.} There exists no such vertex satisfying the above inequality. Then by the proof of Theorem \ref{upper}, the cardinality of $M$ is at most $\max \{k: \sum_{i=1}^k (d_i+1) \leq n\bar{t}\}$. Hence $M$ is the desired monopoly.

\noindent {\bf Case 2.} The algorithm finds a vertex say $v$ such that $|N_{A\cup B}(v)|\leq |N_B(v)|+deg(v)-t(v)$. In this case the proof of Theorem \ref{upper} shows that $M\setminus v$ is still a dynamo. We replace $M$ by $M\setminus v$ and obtain the corresponding sets $A$ and $B$ and go to the next step (i.e. vertex scanning stage). We repeat this procedure and finally obtain a dynamic monopoly $M$ satisfying the condition of the proposition.

\noindent In the following we estimate the running time of the algorithm. Each scanning step takes ${\mathcal{O}}(n^2)$ times. Since at each step one vertex is removed from the dynamic monopoly $M$, the number of steps are at most $n$. It follows
that the total running time is ${\mathcal{O}}(n^3)$. We summarize the algorithm in the following pseudocode form:\\

\begin{algorithmic}
\WHILE{$|N_{A\cup B}(v)| < |N_B(v)|+deg(v)-t(v)+1$}

\STATE $M = M \setminus \{v\}$

\STATE update $A$

\STATE update $B$
\ENDWHILE
\end{algorithmic}
\end{proof}

\noindent By a minimal monopoly $M$ we mean any dynamic monopoly such that no proper subset of $M$ is a dynamic monopoly. We have the following remark from the proof of Theorem \ref{upper}.

\begin{remark}
Let $G$ be a graph with degree sequence $d_1 \leq d_2 \leq \ldots \leq d_n$ in increasing
form. Let also $\tau$ be any threshold assignment for the vertices of $G$ with average threshold $\bar{t}$. Let also $M$ be any minimal $\tau$-dynamic monopoly. Then $$|M|\leq \max \{k: \sum_{i=1}^k (d_i+1) \leq n\bar{t}\}.$$
\end{remark}

\noindent The following is the immediate corollary of Theorem \ref{upper}.

\begin{cor}
Let $G$ be a graph on $n$ vertices and with the minimum degree $\delta$. Then $$
Dyn_{_{\bar{t}=t}}(G)\leq \frac{n\bar{t}}{\delta+1}.$$
\end{cor}

\noindent In the following we determine the exact value of $Dyn_{_{\bar{t}=t}}(K_n)$.

\begin{prop}
$$Dyn_{_{\bar{t}=t}}(K_n)=\lfloor t \rfloor.$$
\end{prop}

\noindent \begin{proof}
Note first that $Dyn_{_{\bar{t}=t}}(K_n)\leq t$ by Theorem \ref{upper}.
In order to prove the converse inequality we consider the following thresholds whose average is $\bar{t}$. Take $n(t-\lfloor t \rfloor)$ vertices of $K_n$ with threshold
equal to $\lfloor t \rfloor +1$ and $n(1-t+\lfloor t \rfloor)$ vertices with threshold $\lfloor t \rfloor$. Note that $nt$ is an integer since $nt=\sum_{v\in K_n} t(v)$. The average of these thresholds is $t$. It is clear that any dynamo for this set of thresholds needs at least $\lfloor t \rfloor$ vertices. This completes the proof.
\end{proof}


\end{document}